\documentclass[reqno]{amsart}

\usepackage{amsmath,amsthm,amssymb,amsfonts,enumerate,graphicx,color,pstricks}
\usepackage[utf8]{inputenc}

\numberwithin{equation}{section} 
\theoremstyle{plain}
\newtheorem{theo+}           {Theorem}      [section]
\newtheorem{prop+}  [theo+]  {Proposition}
\newtheorem{coro+}  [theo+]  {Corollary}
\newtheorem{lemm+}  [theo+]  {Lemma}
\newtheorem{defi+}  [theo+]  {Definition}
\newtheorem{conj+}  [theo+]  {Conjecture}

\theoremstyle{definition}
\newtheorem{rema+}  [theo+]  {Remark}
\newtheorem{prob+}  [theo+]  {Problem}
\newtheorem{exam+}  [theo+]  {Example}

\newenvironment{proposition}{\begin{prop+}}{\end{prop+}}
\newenvironment{corollary}{\begin{coro+}}{\end{coro+}}

\renewcommand{\Re}{\operatorname{Re}}

\begin{document}

\baselineskip 18pt
\larger[2]
\title
[String theory amplitudes and partial fractions] 
{String theory amplitudes and partial fractions}
\author{Hjalmar Rosengren}
\address
{Department of Mathematical Sciences
\\ Chalmers University of Technology and University of Gothenburg\\SE-412~96 G\"oteborg, Sweden}
\email{hjalmar@chalmers.se}
\urladdr{http://www.math.chalmers.se/{\textasciitilde}hjalmar}

\begin{abstract}
We give rigorous proofs and generalizations of partial fraction expansions for string amplitudes that were recently discovered by Saha and Sinha.
\end{abstract}

\maketitle

\section{Introduction}  
Veneziano \cite{v}  interpreted the function
$$B(x_1,x_2)+B(x_1,x_3)+B(x_2,x_3) $$
as a scattering amplitude for strongly interacting mesons. Here, 
$$B(x_1,x_2)=\frac{\Gamma(x_1)\Gamma(x_2)}{\Gamma(x_1+x_2)} $$
is the beta function. Soon afterwards, Virasoro \cite{vi}  and Shapiro \cite{sh} found other combinations of beta functions with similar properties, such as
\begin{equation}\label{via} \frac{\Gamma(x_1)\Gamma(x_2)\Gamma(x_3)}{\Gamma(1-x_1)\Gamma(1-x_2)\Gamma(1-x_3)}, \end{equation}
where $x_1+x_2+x_3=1$. 
These findings were crucial for the early development of string theory, where the Veneziano and Virasoro--Shapiro amplitudes are associated with open and closed bosonic strings, respectively. Related expressions appear in super\-symmetric string theory; for instance, tree level scattering in type II theories is related to \eqref{via}
with $x_1+x_2+x_3=0$ \cite[Eq.\ (7.4.56)]{gsw}.

As a function of $x_1$, the beta function has poles at the 
 non-positive integers and can be written as an infinite partial fraction
\begin{equation}\label{gapf}B(x_1,x_2)=\sum_{k=0}^\infty\frac{(1-x_2)_k}{k!}\cdot\frac 1{x_1+k},\qquad \Re(x_2)>0.
\end{equation}
The identity \eqref{gapf}  is a special case of  Gauss' summation formula for the hypergeometric function ${}_2F_1$ \cite[Thm.\ 2.2.2]{aar}. Here we use the standard notation
$$(a)_k=\frac{\Gamma(a+k)}{\Gamma(a)}=\begin{cases}a(a+1)\dotsm (a+k-1), & k\geq 0,\\
(a-1)^{-1}(a-2)^{-1}\dotsm(a+k)^{-1}, & k<0.\end{cases} $$

An unattractive feature of \eqref{gapf} is that the variables $x_1$ and $x_2$ play a different role on the right-hand side.  In
 \cite{ss}, Saha and Sinha used physics arguments to obtain more symmetric expansions of string amplitudes. To give an example, the case
  $\alpha=\beta=p=0$ of \cite[Eq.\ (4)]{ss}
 can be written
 \begin{equation}\label{ssei}B(x_1,x_2)=\sum_{k=0}^\infty\frac 1{k!}\left(1-\lambda+\frac{(\lambda-x_1)(\lambda-x_2)}{\lambda+k}\right)_k\left(\frac 1{x_1+k}+\frac 1{x_2+k}-\frac 1{\lambda+k}\right). \end{equation}
 Here, $\lambda$ is a free parameter, subject only to the convergence condition $\Re(\lambda)>0$. 
 If $\lambda=x_2$, the symmetry between $x_1$ and $x_2$ is broken and we recover \eqref{gapf}.
Another interesting case is the limit $\lambda\rightarrow\infty.$
Since
$$1-\lambda+\frac{(\lambda-x_1)(\lambda-x_2)}{\lambda+k}=1-k-x_1-x_2+\frac{(x_1+k)(x_2+k)}{\lambda+k}, $$
this limit is
$$B(x_1,x_2)=\sum_{k=0}^\infty\frac {(-1)^k(x_1+x_2)_k}{k!}\left(\frac 1{x_1+k}+\frac 1{x_2+k}\right). $$
This is a special case of a known ${}_4F_3$ summation \cite[Cor.\ 3.5.3]{aar}.

As noted by Saha and Sinha, special cases of their identities give intriguing new 
formulas  for $\pi$ and other mathematical constants. For instance, 
  the case $x_1=x_2=1/2$ of \eqref{ssei} is 
$$\pi=\sum_{k=0}^\infty\frac 1{k!}\left(1-\lambda+\frac{(\lambda-1/2)^2}{\lambda+k}\right)_k\left(\frac 4{2k+1}-\frac 1{\lambda+k}\right). $$
For  $\lambda=\infty$, this is the well-known Madhava series
$$\frac\pi 4=\sum_{k=0}^\infty\frac{(-1)^k}{2k+1},  $$
whereas the the case $\lambda=1/2$ is the identity
$$ \frac\pi2=\sum_{k=0}^\infty\frac{(2k-1)!!}{k!2^k(2k+1)},$$
which can be obtained by letting $x=1$ in the Taylor series for $\arcsin(x)$.

Saha and Sinha also gave several expansions related to closed strings. 
For instance, the case $\alpha=1/3$, $\beta=2/3$ of \cite[Eq.\ (A3)]{ss}
gives the following expansion of the Virasoro--Shapiro amplitude \eqref{via}:
\begin{equation}\label{ssce}
 \frac{\Gamma(x_1)\Gamma(x_2)\Gamma(x_3)}{\Gamma(1-x_1)\Gamma(1-x_2)\Gamma(1-x_3)}
 =\sum_{k=0}^\infty\frac{(\eta_k)_k^2}{(k!)^2}
 \left(\sum_{j=1}^3\frac 1{x_j+k}-\frac {1}{\frac 13+k}\right).
\end{equation}
Here,  $x_1+x_2+x_3=1$ and
$$\eta_k=\frac{1-k}2+\sqrt{
\frac{\left(\frac 13+k\right)^2}4+\frac{(x_1-\frac 13)(x_2-\frac 13)(x_3-\frac 13)}{\frac 13+k}}. $$
Since
\begin{align}\nonumber\left(\frac{1-k}2+\sqrt a\right)_k^2&= (-1)^k\left(\frac{1-k}2+\sqrt a\right)_k\left(\frac{1-k}2-\sqrt a\right)_k\\
\label{psc}&=\prod_{j=0}^{k-1}\left(a-\left(\frac{1-k}2+j\right)^2\right),\end{align}
the terms in \eqref{ssce} are rational functions of the parameters.

The purpose of the present work is to give independent and rigorous proofs of the mathematical results of \cite{ss}. We also give  generalizations, which we hope may have some interest for physicists. For instance, when  $x_1+x_2+x_3=s$ is an integer and
$\Re(\lambda+1)>0$, we show that 
\begin{multline}\label{lvia}
 \frac{\Gamma(x_1)\Gamma(x_2)\Gamma(x_3)}{\Gamma(1-x_1)\Gamma(1-x_2)\Gamma(1-x_3)}\\
 =\sum_{k=0}^\infty\frac{(-1)^{s+1}(\xi_k
 )_{k+s-1}^2}{(k!)^2}
 \left(\sum_{j=1}^3\frac 1{x_j+k}-\frac {1}{\lambda+k}\right),
\end{multline}
where
$$\xi_k=\frac{2-s-k}2+\sqrt{\left(\frac{s+k-2\lambda}2\right)^2+\frac{(x_1-\lambda)(x_2-\lambda)(x_3-\lambda)}{\lambda+k}}. $$
As mentioned above, the cases $s=0$ and $s=1$ are of special relevance in physics. In \cite{ss}, \eqref{lvia} is only obtained for the fixed value $\lambda=s/3$. In particular, the case $s=1$ and $\lambda=1/3$ is \eqref{ssce}.

Our method is to obtain infinite partial fractions as limit cases of finite partial fractions for symmetric rational functions.
One source of inspiration was Schlosser's work \cite[\S 7]{s}, which is the only place where I have seen series remotely similar to \eqref{ssei}. The matrix inversions used by Schlosser are in fact closely related to partial fractions \cite{rs}, so there may be more direct connections that remain to be explored. 

{\bf Acknowledgements:} I thank Aninda Sinha for illuminating correspondence. I thank the community at MathOverflow, especially Martin Nicholson (user Nemo) for a crucial comment. This research is supported by the Swedish Science Research council, project 2020-04221.

\section{Partial fraction expansions for symmetric rational functions}

As a warm-up, we consider the identity \eqref{gapf}, where we write $x_1=x$, $x_2=a$. One way to derive it is to start from the finite truncation
 $$\frac{\Gamma(x)}{\Gamma(a+x)}\cdot\frac{\Gamma(a+x+n)}{\Gamma(x+n+1)}=\frac{(a+x)_n}{(x)_{n+1}}.  $$
It has the partial fraction expansion
\begin{equation}\label{bpf}\frac{(a+x)_n}{(x)_{n+1}}= \sum_{k=0}^n \frac{B_k}{x+k}, \end{equation}
where
$$B_k=\lim_{x\rightarrow -k}(x+k) \frac{(a+x)_n}{(x)_{n+1}}=
\frac{(-1)^k(a-k)_n}{k!(n-k)!}.
$$
The identity \eqref{bpf} can be recognized as a special case of the Pfaff--Saalsch\"utz ${}_3F_2$ summation \cite[Thm.\ 2.2.6]{aar}.

We now let $n\rightarrow\infty$ in \eqref{bpf}.
Using that
\begin{equation}\label{wsf}\frac{\Gamma(a+n)}{\Gamma(b+n)}\sim n^{a-b},\qquad n\rightarrow\infty, \end{equation}
it is straight-forward to formally obtain \eqref{gapf} in the limit. To make this rigorous one can apply Tannery's theorem 
if $\Re(x_2)>1$, and then extend the identity analytically to the larger domain $\Re(x_2)>0$. This is discussed for more general series in the Appendix.

 To obtain an analogous proof of \eqref{ssei}, we start from
$$ \frac{\Gamma(x_1)\Gamma(x_2)}{\Gamma(x_1+x_2)}\cdot  \frac{\Gamma(x_1+x_2+n)}{\Gamma(x_1+n+1)\Gamma(x_2+n+1)}
=\frac{(x_1+x_2)_n}{(x_1)_{n+1}(x_2)_{n+1}}.
$$
We then apply a partial fraction expansion for symmetric rational functions in two variables. For \eqref{ssei}, we need to work with three variables. In the remainder of this section, we
discuss the general case of $r$ variables.

 Recall that the ring of symmetric polynomials is generated by the
 elementary symmetric polynomials, which may be defined by the generating function
 \begin{equation}\label{egf}\prod_{j=1}^r(t+x_j)=\sum_{m=0}^r t^{r-m} e_m(\mathbf x).\end{equation}
Thus, to compute the value of a symmetric polynomial at a point $\mathbf x=(x_1,\dots,x_r)$, it suffices to know
$e_m(\mathbf x)$ for $1\leq m\leq r$.

\begin{proposition}\label{gpfp}
Let $P$ be a symmetric polynomial of $r$ variables, which is of degree at most $n+1$ in each variable. Let $Q_0,\dots,Q_n$ be symmetric polynomials that are of degree exactly $1$ in each variable. Then,
\begin{equation}\label{gpf}\frac{P(\mathbf x)}{\prod_{j=0}^nQ_j(\mathbf x)}-\frac{P(\mathbf y)}{\prod_{j=0}^nQ_j(\mathbf y)}
=\sum_{k=0}^n\frac{P(\mathbf t_k)}{\prod_{j=0,j\neq k}^nQ_j(\mathbf t_k)}\left(\frac 1{Q_k(\mathbf x)}-\frac 1{Q_k(\mathbf y)}\right),
 \end{equation}
 where $\mathbf t_k$ denotes the vector defined up to permutations by
 \begin{equation}\label{emtd}e_m(\mathbf t_k)=\frac{Q_k(\mathbf y)e_m(\mathbf x)-Q_k(\mathbf x)e_m(\mathbf y)}{Q_k(\mathbf y)-Q_k(\mathbf x)}. \end{equation}
\end{proposition}

By \eqref{egf}, we can equivalently write \eqref{emtd} as
\begin{equation}\label{tad} (Q_k(\mathbf y)-Q_k(\mathbf x))\prod_{j=1}^r(u-t_j)=Q_k(\mathbf y) \prod_{j=1}^r(u-x_j)-Q_k(\mathbf x)\prod_{j=1}^r(u-y_j),\end{equation}
where $\mathbf t_k=(t_1,\dots,t_r)$. That is,  the components of $\mathbf t_k$ 
 are  the solutions to
\begin{equation}\label{tae}Q_k(\mathbf x)\prod_{j=1}^r(u-y_j)=Q_k(\mathbf y) \prod_{j=1}^r(u-x_j),
\end{equation}
considered as an equation in $u$. 

\begin{proof}
Since we are proving a rational function identity, we may assume that all parameters are generic.
 Define
$\lambda_m$ and $\mu_m$ by
\begin{subequations}\label{exyg}
\begin{align}\label{exy}e_m(\mathbf x)&=\lambda_m e_1(\mathbf x)+\mu_m,\\
 \label{ey}e_m(\mathbf y)&=\lambda_m e_1(\mathbf y)+\mu_m, \end{align}
\end{subequations}
that is,
$$\lambda_m=\frac{e_m(\mathbf x)-e_m(\mathbf y)}{e_1(\mathbf x)-e_1(\mathbf y)},\qquad 
\mu_m=\frac{e_1(\mathbf x)e_m(\mathbf y)-e_m(\mathbf x)e_1(\mathbf y)}{e_1(\mathbf x)-e_1(\mathbf y)}. 
$$
  Considering $\lambda_m$ and $\mu_m$ as constants, we can view any symmetric polynomial in $\mathbf x$ as a function of the single variable $e_1(\mathbf x)$.
  
  For fixed $k$, let us write $Q_k=Q$ and $\mathbf t_k=\mathbf t$. 
  Note that
  \begin{equation}\label{qae}Q(\mathbf x)=\sum_{m=0}^r A_m  e_m(\mathbf x)\end{equation}
  for some coefficients $A_m$.
  Multiplying \eqref{emtd} with $A_m$ and summing over $m$ gives
  \begin{equation}\label{qtv}Q(\mathbf t)=0.\end{equation}
  It is also clear from \eqref{emtd} and \eqref{exyg} that
  \begin{equation}\label{emt}e_m(\mathbf t)=\lambda_m e_1(\mathbf t)+\mu_m. \end{equation}
  The equations \eqref{qtv} and \eqref{emt} determine $\mathbf t$ uniquely. Hence,
   $\mathbf t$ depends on $\mathbf x$ and $\mathbf y$ only through the parameters $\lambda_m$ and $\mu_m$. 
  To make this explicit, we write
$$0=\sum_{m=0}^r A_{m}e_m(\mathbf t)=\sum_{m=0}^r A_{m}(\lambda_m e_1(\mathbf t)+\mu_m).
$$
Solving this for $e_1(\mathbf t)$ and using again \eqref{emt} gives
\begin{equation}\label{emtb}e_m(\mathbf t)=-\lambda_m\frac{\sum_{j=0}^r A_j\mu_j}{\sum_{j=0}^r A_j\lambda_j}+\mu_m. 
\end{equation}

We are now reduced to a  partial fraction expansion in one variable. Let $L$ and $R$ denote the two sides of \eqref{gpf} and
 $S=\prod_{j=0}^n Q_j$. Using \eqref{exy}, we can write $SL$ as a polynomial in  
$e_1(\mathbf x)$ of  degree at most $n+1$. 
By \eqref{emtb}, the coefficient $P(\mathbf t_k)/{\prod_{j\neq k}Q_j(\mathbf t_k)}$  can be treated as a constant. Hence, $SR$ is also
 a polynomial in  
$e_1(\mathbf x)$ of  degree at most $n+1$.

To complete the proof, it suffices to check that $SL$ and $SR$ agree for $n+2$ independent values of $\mathbf x$, subject to  \eqref{exy}. By \eqref{ey} and \eqref{emt}, we can take these as $\mathbf t_0,\dots, \mathbf  t_n,\mathbf y$. By 
\eqref{qtv}, if $\mathbf x=\mathbf t_k$,  only the
first term on the left and the 
$k$-th term on the right contribute to the value of $SR$. It is then straight-forward to check that $SL=SR$. The identity  for $\mathbf x=\mathbf y$ is obvious.
This completes the proof.
\end{proof}

Let us now consider the limit case of Proposition \ref{gpfp} 
 when 
 $y_1,\dots,y_l$ are fixed and $y_{l+1},\dots,y_r\rightarrow\infty$, where $0\leq l\leq r-1$.
  If 
 $P$ is of degree at most $n$ in each variable, then 
 all terms in \eqref{gpf} that explicitly depend on $\mathbf y$  tend to zero. 
 In particular, the left-hand side is independent of $y_1,\dots,y_{l}$, so we can think of these as free parameters.  
 If $Q$ is as in \eqref{qae}, we  write
 $$\hat Q(y_1,\dots,y_{l})=\lim_{y_{l+1},\dots,y_r\rightarrow\infty}\frac {Q(y_1,\dots,y_r)}{y_{l+1}\dotsm y_r}
 =\sum_{m=0}^l A_{m+r-l} \,e_{m}(y_1,\dots,y_l). $$
 Then,  the equation \eqref{emtd} degenerates to
 \begin{equation}\label{emti}e_m(\mathbf t_k)=e_m(\mathbf x)-\frac{Q_k(\mathbf x)}{\hat Q_k(\mathbf y)}\,e_{m+l-r}(\mathbf y),
 \end{equation}
 where we should interpret $e_{k}(\mathbf y)$ as $0$ for $k<0$.

 \begin{corollary}
Let $P$ and $Q_j$ be as in Proposition \ref{gpfp}, but assume that $P$ is of degree at most $n$ in each variable. Let $y_1,\dots,y_{l}$ be generic  scalars and let $\mathbf t_k$ be defined up to permutations by \eqref{emti}.
Then,
\begin{equation}\label{gpfi}\frac{P(\mathbf x)}{\prod_{j=0}^nQ_j(\mathbf x)}
=\sum_{k=0}^n\frac{P(\mathbf t_k)}{\prod_{j=0,j\neq k}^nQ_j(\mathbf t_k)}\cdot\frac 1{Q_k(\mathbf x)}.
 \end{equation}
\end{corollary}

Returning to Proposition \ref{gpfp}, suppose that $Q_k(\mathbf x)=\prod_{m=1}^r(x_m-b_k)$
for some scalars $b_k$. Then, \eqref{tae} reduces to
$$\prod_{j=1}^r(b_k-x_j)(u-y_j)= \prod_{j=1}^r(b_k-y_j)(u-x_j),$$
which has one obvious solution $u=b_k$. We write the full vector of solutions as $\mathbf t_k=(b_k,b_k',\dots, b_k^{(r-1)})$. 
Differentiating \eqref{tad} in $u$ 
and then substituting $u=b$ gives after simplification
$$\frac 1{Q(\mathbf x)}-\frac 1{Q(\mathbf y)}=\frac 1{\prod_{j=1}^{r-1}(b^{(j)}-b)}\left(\sum_{j=1}^r\frac 1{x_j-b}-\sum_{j=1}^r\frac 1{y_j-b}\right). $$
Thus, \eqref{gpf} can be written
\begin{multline}\label{gfd}\frac{P(\mathbf x)}{\prod_{j=0}^n\prod_{m=1}^r(x_m-b_j)}-\frac{P(\mathbf y)}{\prod_{j=0}^n\prod_{m=1}^r(y_m-b_j)}\\
=\sum_{k=0}^n\frac{P(b_k,b_k',\dots,b_k^{(r-1)})}{\prod_{j=0,j\neq k}^n(b_k-b_j)\prod_{j=0}^n\prod_{m=1}^{r-1}(b_k^{(m)}-b_j) }\left(\sum_{j=1}^r\frac 1{x_j-b_k}-\sum_{j=1}^r\frac 1{y_j-b_k}
\right).
 \end{multline}
 Taking as before the limit  $y_{l+1},\dots,y_r\rightarrow\infty$, when
  $P$ has degree at most $n$ in each variable, gives 
\begin{multline}\label{fipf}\frac{P(\mathbf x)}{\prod_{j=0}^n\prod_{m=1}^r(x_m-b_j)}
=\sum_{k=0}^n\frac{P(b_k,b_k',\dots,b_k^{(r-1)})}{\prod_{j=0,j\neq k}^n(b_k-b_j)\prod_{j=0}^n\prod_{m=1}^{r-1}(b_k^{(m)}-b_j) }\\
\times\left(\sum_{j=1}^r\frac 1{x_j-b_k}-
\sum_{j=1}^{l}\frac 1{y_j-b_k}
\right).
 \end{multline}
 Here, $\mathbf t_k=(b_k,b_k',\dots,b_k^{(r-1)})$ is given by \eqref{emti}, which in this case reads
 \begin{equation}\label{emtf}e_m(\mathbf t_k)=e_m(\mathbf x)-\frac{\prod_{j=1}^r(x_j-b_k)}{\prod_{j=1}^{l}(y_j-b_k)}\,
 e_{m+l-r}(\mathbf y) \end{equation}
 or, equivalently,
 \begin{equation}\label{emtgf}\prod_{j=1}^r(b_k^{(j)}-u)=\prod_{j=1}^r(x_j-u)-
 \prod_{j=1}^{l}\frac{y_j-u}{y_j-b_k}\prod_{j=1}^r(x_j-b_k). \end{equation}

We will in fact only need the  cases $(r,l)=(2,1)$ and $(r,l)=(3,1)$ of the results above. However, we find it instructive and potentially useful to state them in greater generality.

\section{Expansions of open string amplitudes}
\label{oss}

To prove \eqref{ssei} and some related results,  we start from the case $r=2$, $l=1$  of \eqref{fipf}. Writing $y_1=\lambda$, it has the form
\begin{multline}\label{twpf}
\frac{P(x_1,x_2)}{\prod_{j=0}^n(x_1-b_j)(x_2-b_j)}\\
=\sum_{k=0}^n\frac{P(b_k,b_k')}{\prod_{j=0,\,j\neq k}^n(b_k-b_j)\prod_{j=0}^n( b_k'-b_j)}\left(\frac 1{x_1-b_k}+\frac 1{x_2-b_k}-\frac 1{\lambda-b_k}\right).
 \end{multline}
By the case  $u=\lambda$ of \eqref{emtgf},
$$ (b-\lambda)(b'-\lambda)=(x_1-\lambda)(x_2-\lambda),$$
which gives 
$$b'
=\lambda-\frac{(\lambda-x_1)(\lambda-x_2)}{\lambda-b}.$$

We now specialize to the case
   $b_k=-k$ and $P(x_1,x_2)=(x_1+x_2+a)_n$.  
   After a straight-forward computation, we arrive at the  rational function identity
  \begin{multline}\label{cbi}\frac{(x_1+x_2+a)_{n}}{(x_1)_{n+1}(x_2)_{n+1}}\\
=\sum_{k=0}^n \frac{(-1)^k}{k!(n-k)!}\frac{\left(\lambda-\frac{(\lambda-x_1)(\lambda-x_2)}{\lambda+k}+a-k\right)_n}{\left(\lambda-\frac{(\lambda-x_1)(\lambda-x_2)}{\lambda+k}\right)_{n+1}}\left(\frac 1{x_1+k}+\frac 1{x_2+k}-\frac 1{\lambda+k}\right).\end{multline}
  
Next, we multiply both sides of \eqref{cbi} by $\Gamma(n+2-a)$ and let
 $n\rightarrow\infty$. Formally using \eqref{wsf}, it is straight-forward to 
  obtain the following result. However, this involves an interchange of limit and summation that needs to be justified. We provide the details in the Appendix.

\begin{corollary}\label{cbl}
 Assuming $\operatorname{Re}(a+\lambda)>0$,
  \begin{multline}\label{cdi}\frac{\Gamma(x_1)\Gamma(x_2)}{\Gamma(x_1+x_2+a)}\\
=\sum_{k=0}^\infty \frac{(-1)^k}{k!}\cdot\frac{\Gamma\left(\lambda-\frac{(\lambda-x_1)(\lambda-x_2)}{\lambda+k}
\right)}{\Gamma\left(\lambda-\frac{(\lambda-x_1)(\lambda-x_2)}{\lambda+k}+a-k\right)}\left(\frac 1{x_1+k}+\frac 1{x_2+k}-\frac 1{\lambda+k}\right).\end{multline}
  \end{corollary}

When 
$a$ is an integer,
\eqref{cdi} can be written
\begin{multline*}\frac{\Gamma(x_1)\Gamma(x_2)}{\Gamma(x_1+x_2+a)}\\
=\sum_{k=0}^\infty \frac{(-1)^a\left(1-\lambda+\frac{(\lambda-x_1)(\lambda-x_2)}{\lambda+k}\right)_{k-a}}{k!}\left(\frac 1{x_1+k}+\frac 1{x_2+k}-\frac 1{\lambda+k}\right).\end{multline*}
 This  is  \cite[Eq.\ (4)]{ss}, although it is only stated there under the additional assumption  $a\leq 1$. 
   The case $a=0$ is the example \eqref{ssei} given in the introduction.

If one is interested in formulas for $\pi$ one can  specialize $x_1=1/2+m$, $x_2=1/2+n$  and use
$$\Gamma\left(\frac 12+m\right)=(1/2)_m\,\Gamma(1/2)=\frac{(2m-1)!!}{2^m}\,\sqrt{\pi}. $$
This gives a family of identities
\begin{multline*}\pi=\frac{(-1)^a(m+n+a)!2^{m+n}}{(2m-1)!!(2n-1)!!}\sum_{k=0}^\infty \frac{1}{k!}
\left(1-\lambda+\frac{(\lambda-\frac 12-m)(\lambda-\frac 12-n)}{\lambda+k}\right)_{k-a}\\
\times
\left(\frac 2{2m+2k+1}+\frac 2{2n+2k+1}-\frac 1{\lambda+k}\right),
\end{multline*}
parametrized by three integers $m$, $n$ and $a$  with $m\geq 0$, $n\geq 0$ and $m+n+a\geq 0$ and a continuous parameter $\lambda$  with  $\Re(\lambda)>-a$.

 \section{Asymmetric expansions of closed string amplitudes}
 \label{css}

Saha and Sinha gave two types of expansions for amplitudes related to closed strings. We start with the first type
 \cite[Eq.\ (A2)]{ss}, which is symmetric in the variables 
 $x_1$ and $x_2$ but where $x_3$ plays a different role. Taking the Virasoro--Shapiro amplitude \eqref{via} as an example, we will treat $x_1$ and $x_2$ as variables, but $x_3$ as determined from the relation $x_1+x_2+x_3=1$. 
The left-hand side then has some poles of the form $x_1=x_2=-k$ (for $k$ a non-negative integer) and others of the form $x_1+x_2=k+1$. We will therefore start from symmetric rational functions of the form
$$\frac{P(x_1,x_2)}{\prod_{j=0}^n(x_1-b_j)(x_2-b_j) (c_j-x_1-x_2)}. $$

We will apply \eqref{gpfi} with $r=2$, $l=1$ and write
 $y_1=\lambda$.
  We then have $n+1$ terms 
 that can be written as in \eqref{twpf}, with the additional denominator factors
 $ \prod_{j=0}^n (c_j-b_k-b_k')$. Let us write the vector
 $\mathbf t$ corresponding to a factor $Q=c-x_1-x_2$ as $(c^+,c^-)$. The equations \eqref{emti} reduce to
 $$c^++c^-=c,\qquad c^+c^-=x_1x_2+(c-x_1-x_2)\lambda. $$
  Hence,
 $$c^\pm=\frac{c^++c^-}{2}\pm\sqrt{\frac{(c^++c^-)^2-4c^+c^-}{4}}=\frac{c}{2}\pm\sqrt{\frac{c^2}4+\lambda(x_1+x_2- c)-x_1x_2}. $$
This leads to the identity
\begin{multline*}
\frac{P(x_1,x_2)}{\prod_{j=0}^n(x_1-b_j)(x_2-b_j)(c_j-x_1-x_2)}\\
=\sum_{k=0}^n\frac{ P(b_k,b_k')}{\prod_{j=0,\,j\neq k}^n (b_k-b_j)\prod_{j=0}^n(b_k'-b_j)(c_j-b_k-b_k')}\\
\times\left(\frac 1{x_1-b_k}+\frac 1{x_2-b_k}-\frac 1{\lambda-b_k}\right)\\
+\sum_{k=0}^n\frac{ P(c_k^+,c_k^-)}{\prod_{j=0}^n (c_k^+-b_j)(c_k^--b_j)\prod_{j=0,\,j\neq k}^n(c_j-c_k)}\cdot\frac 1{c_k-x_1-x_2},
 \end{multline*}
which holds for $P$ a symmetric polynomial of degree at most $2n+1$ in each variable.

Specializing
$b_k=-k$, $c_k=s+k$ and $P(x_1,x_2)=(t-x_1)_n(t-x_2)_n(u+x_1+x_2)_n$ gives after simplification
\begin{multline}\label{affp}
\frac{(t-x_1)_n(t-x_2)_n(u+x_1+x_2)_n}{(x_1)_{n+1}(x_2)_{n+1}(s-x_1-x_2)_{n+1}}\\
\begin{split}&=\sum_{k=0}^n\frac{(-1)^k}{k!(n-k)!}\frac{(t+k)_n\left(t-
\lambda+\frac{(\lambda-x_1)(\lambda-x_2)}{\lambda+k} 
\right)_n\left(u+\lambda-\frac{(\lambda-x_1)(\lambda-x_2)}{\lambda+k} -k\right)_n}{\left(\lambda-\frac{(\lambda-x_1)(\lambda-x_2)}{\lambda+k} \right)_{n+1}\left(s-\lambda+\frac{(\lambda-x_1)(\lambda-x_2)}{\lambda+k} +k\right)_{n+1}}\\
&\quad\quad\times\left(\frac 1{x_1+k}+\frac 1{x_2+k}-\frac 1{\lambda+k}\right)\\
&\quad+\sum_{k=0}^n\frac{(-1)^k}{k!(n-k)!}\frac{(t-c_k^+)_n(t-c_k^-)_n(u+s+k)_n}{(c_k^+)_{n+1}(c_k^-)_{n+1}}\cdot\frac 1{s+k-x_1-x_2}.
\end{split}\end{multline}
If we formally let $n\rightarrow\infty$, we obtain the following result.
To make this rigorous, we again need some estimates that are given in 
 the Appendix.

\begin{corollary}\label{acc}
Assuming $\Re(2\lambda-s+t+u)>0$, we have
\begin{multline}\label{clf}
\frac{\Gamma(x_1)\Gamma(x_2)\Gamma(s-x_1-x_2)}{\Gamma(t-x_1)\Gamma(t-x_2)\Gamma(u+x_1+x_2)}
\\
\begin{split}&=\sum_{k=0}^\infty\frac{(-1)^k}{k!}\frac{\Gamma\left(\lambda-\frac{(\lambda-x_1)(\lambda-x_2)}{\lambda+k}\right)\Gamma\left(s-\lambda+\frac{(\lambda-x_1)(\lambda-x_2)}{\lambda+k}+k\right)}{\Gamma(t+k)\Gamma\left(t-\lambda+\frac{(\lambda-x_1)(\lambda-x_2)}{\lambda+k}\right)\Gamma\left(u+\lambda-\frac{(\lambda-x_1)(\lambda-x_2)}{\lambda+k}-k\right)}\\
&\quad\quad\times\left(\frac 1{x_1+k}+\frac 1{x_2+k}-\frac 1{\lambda+k}\right)\\
&\quad+\sum_{k=0}^\infty\frac{(-1)^k}{k!}\frac{\Gamma(c_k^+)\Gamma(c_k^-)}{\Gamma(t-c_k^+)\Gamma(t-c_k^-)\Gamma(u+s+k)}\cdot\frac 1{s+k-x_1-x_2},
\end{split}\end{multline}
where
$$c_k^\pm=\frac{s+k}{2}\pm\sqrt{\frac{(s+k)^2}4+\lambda(x_1+x_2-s-k)-x_1x_2}. $$
\end{corollary}

As an example, let $s=t=1$ and $u=0$ in \eqref{clf} and reintroduce the variable  $x_3=1-x_1-x_2$.
After simplification,
this gives the following asymmetric expansion of the Virasoro--Shapiro amplitude \eqref{via}:
\begin{multline*}
\frac{\Gamma(x_1)\Gamma(x_2)\Gamma(x_3)}{\Gamma(1-x_1)\Gamma(1-x_2)\Gamma(1-x_3)}
\\
=\sum_{k=0}^\infty\frac{\left(1-\lambda+\frac{(\lambda-x_1)(\lambda-x_2)}{\lambda+k}\right)_k^2}{(k!)^2}\left(\frac 1{x_1+k}+\frac 1{x_2+k}-\frac 1{\lambda+k}\right)\\
+\sum_{k=0}^\infty\frac{\left(\frac{1-k}{2}+\sqrt{\frac{(1+k)^2}4-\lambda(x_3+k)-x_1x_2}\right)_k^2}{(k!)^2}\cdot\frac 1{x_3+k},
\end{multline*}
where $x_1+x_2+x_3=1$ and $\Re(\lambda)>0$. By \eqref{psc}, the terms are rational functions of the parameters.
If we instead let $s=0$ and $t=u=1$  we recover \cite[Eq.\ (A2)]{ss}.

\section{Symmetric expansions of closed string amplitudes}
 \label{scs}

To obtain  symmetric expansions for closed string amplitudes, such as \eqref{ssce}, we start from 
$$\frac{(u-x_1)_n(u-x_2)_n(u-x_3)_n}{(x_1)_{n+1}(x_2)_{n+1}(x_3)_{n+1}}
$$
and apply \eqref{fipf} with $r=3$ and $l=1$. 
We will write $y_1=\lambda$ and $t^+$, $t^-$ instead of $t'$, $t''$. It follows from \eqref{emtf} that 
$$t+t^++t^-=x_1+x_2+x_3 $$
and from the case $u=\lambda$ of \eqref{emtgf} that
$$(t-\lambda)(t^+-\lambda)(t^--\lambda)=(x_1-\lambda)(x_2-\lambda)(x_3-\lambda). $$
This leads to the explicit formula
\begin{align*}t^\pm&=\frac{t^++t^-}{2}\pm\sqrt{\left(\frac{t^++t^--2\lambda}{2}\right)^2-
(t^+-\lambda)(t^--\lambda)}\\
&=\frac{e_1(\mathbf x)-t}{2}\pm\sqrt{\frac {({e_1(\mathbf x)-t-2\lambda})^2}4-
\frac{\prod_{j=1}^3(x_j-\lambda)}{t-\lambda}}.\end{align*}
The resulting special case of \eqref{fipf} is
\begin{multline}\label{yl}\frac{(u-x_1)_n(u-x_2)_n(u-x_3)_n}{(x_1)_{n+1}(x_2)_{n+1}(x_3)_{n+1}}\\
=\sum_{k=0}^n\frac{(-1)^k}{k!(n-k)!}\frac{(u+k)_n(u-\eta_k^+)_n(u-\eta_k^-)_n}{(\eta_k^+)_{n+1}(\eta_k^-)_{n+1}}\left(\sum_{j=1}^3\frac 1{x_j+k}-
\frac 1{\lambda+k}\right),
\end{multline}
where we write $\eta_k^\pm=(-k)^{\pm}$. 
If we let $n\rightarrow\infty$, we obtain the following result.
As before, we provide some details in the Appendix.

\begin{corollary}\label{gcac}
If $\Re(\lambda+u-e_1(\mathbf x))>0$ then
  \begin{multline}\label{tvd}
 \frac{\Gamma(x_1)\Gamma(x_2)\Gamma(x_3)}{\Gamma(u-x_1)\Gamma(u-x_2)\Gamma(u-x_3)}\\
 =\sum_{k=0}^\infty\frac{(-1)^k}{k!}\frac{\Gamma(\eta_k^+)\Gamma(\eta_k^-)}{\Gamma(u+k)\Gamma(u-\eta_k^+)\Gamma(u-\eta_k^-)}
 \left(\sum_{j=1}^3\frac 1{x_j+k}-\frac 1{\lambda+k}\right),
\end{multline}
where 
$$\eta_k^\pm=\frac{e_1(\mathbf x)+k}{2}\pm\sqrt{\frac {({e_1(\mathbf x)+k-2\lambda})^2}4+
\frac{\prod_{j=1}^3(x_j-\lambda)}{\lambda+k}}.$$
\end{corollary}

 If we assume that $u-e_1(\mathbf x)=a$ is an integer,  \eqref{tvd} can be written
as
\begin{multline}\label{tve}
 \frac{\Gamma(x_1)\Gamma(x_2)\Gamma(x_3)\Gamma(a+s)}{\Gamma(a+s-x_1)\Gamma(a+s-x_2)\Gamma(a+s-x_3)}\\
 =\sum_{k=0}^\infty\frac{(-1)^k}{k!}\frac{(1-\eta_k^+)_{k-a}(1-\eta_k^-)_{k-a}}{(a+s)_k}
 \left(\sum_{j=1}^3\frac 1{x_j+k}-\frac {1}{\lambda+k}\right),
\end{multline}
where $x_1+x_2+x_3=s$.  The special case $\lambda=s/3$ of \eqref{tve} is \cite[Eq.\ (A3)]{ss}. The case  $a+s=1$  is
the identity \eqref{lvia} given in the introduction.

As another example, the case $x_1=x_2=x_3=1/2$, $a=0$ of \eqref{tve} can be written
$$\frac{\pi^2}2=\sum_{k=0}^\infty\frac{(-1)^k(\xi_k^+)_k(\xi_k^-)_k}{k!(3/2)_k}\left(\frac 6{2k+1}-\frac 1{\lambda+k}\right), \qquad \Re(\lambda)>0,$$
where
$$\xi_k^\pm=\frac 14-\frac k2\pm\frac{2k+1}4\sqrt{\frac{k+2-3\lambda}{k+\lambda}}.$$
When $\lambda=1/2$, this reduces to Euler's identity
$$\frac{\pi^2}8=\sum_{k=0}^\infty\frac 1{(2k+1)^2}. $$
The series obtained by formally letting $\lambda\rightarrow\infty$ is divergent.
 
 The reader may ask why we did not use \eqref{fipf} in the more general case $l=2$.  We would then replace \eqref{yl} with 
 \begin{multline}\label{yyi}\frac{(u-x_1)_n(u-x_2)_n(u-x_3)_n}{(x_1)_{n+1}(x_2)_{n+1}(x_3)_{n+1}}\\
=\sum_{k=0}^n\frac{(-1)^k}{k!(n-k)!}\frac{(u+k)_n(u-(-k)')_n(u-(-k)'')_n}{((-k)')_{n+1}((-k)'')_{n+1}}\\
\times\left(\sum_{j=1}^3\frac 1{x_j+k}-
\sum_{j=1}^2\frac 1{y_j+k}\right).
\end{multline}
 The problem 
is that,  when $n\rightarrow\infty$, the left-hand side behaves as
$n^{3u-2e_1(\mathbf x)-3}$ 
and the terms on the right as $n^{3u-2e_1(\mathbf t)-3}$,
where $\mathbf t=(-k,(-k)',(-k)'')$. 
Hence, to take a termwise limit, we need 
$e_1(\mathbf x)=e_1(\mathbf t)$. However, \eqref{emtf} gives
$$e_1(\mathbf t)=e_1(\mathbf x)-\frac{\prod_{j=1}^3(x_j+k)}{\prod_{j=1}^2(y_j+k)}. $$
This forces  $y_1$ or $y_2\rightarrow\infty$, which is precisely the case $l=1$ studied above. We can also note that, in the special case $y_1=x_1$, $y_2=x_2$, 
it follows from \eqref{emtgf} that $\mathbf t=(-k,x_1,x_2)$.
Then,
\eqref{yyi} reduces to
$$ \frac{(u-x_3)_n}{(x_3)_{n+1}}
=\sum_{k=0}^n\frac{(-1)^k(u+k)_n}{k!(n-k)!}\cdot \frac 1{x_3+k},
$$
which is just a restatement of \eqref{bpf} (with $a=1-u-n$).  In hypergeometric notation, the series on the right is a multiple of 
$${}_3F_2\left(\begin{matrix}-n,u+n,x_3\\u,x_3+1\end{matrix};1\right). $$
The asymptotics as $n\rightarrow\infty$ of this type of series is well studied, see e.g.\ \cite[\S 7.4]{l}. It is conceivable that a similar analysis of the series  \eqref{yyi} would lead to some interesting asymptotic expansions of string amplitudes. 
 
\section*{Appendix. Limit transitions}
\label{app}

We have obtained our main results, Corollary \ref{cbl}, Corollary \ref{acc} and Corollary~\ref{gcac}, by a  limit transition from finite to infinite series. It requires some work to make this rigorous. In the present Appendix, we explain this in detail for
Corollary \ref{cbl} and then briefly discuss the necessary modifications for the other two results.

We start from Stirling's formula
$$\Gamma(z)\sim \sqrt{2\pi}\,z^{z-\frac 12}e^{-z}, $$
which holds if $|z|\rightarrow\infty$ with
$\arg(z)<\pi-\delta$ for any fixed $\delta>0$. It implies that \eqref{wsf} can be strengthened to
\begin{equation}\label{gas}\frac{\Gamma(z+a(z))}{\Gamma(z+b(z))}\sim z^{a(z)-b(z)}, \end{equation}
which holds in the same sense, assuming that $a(z)$ and $b(z)$ are bounded.

As we have already noted, if we multiply both sides of  \eqref{cbi} with $\Gamma(n+2-a)$ and then let $n\rightarrow\infty$ 
 we obtain \eqref{cdi}, provided that we are allowed to interchange  limit and summation on the right-hand side. To justify this, we will use  Tannery's theorem, which states that if 
 $|a_{kn}|<M_k$, where $M_k$ is independent of $n$ and $\sum_{k=0}^\infty M_k<\infty$, then
$$\lim_{n\rightarrow\infty}\sum_{k=0}^\infty a_{kn}=\sum_{k=0}^\infty \lim_{n\rightarrow\infty}a_{kn}. $$
We need to find such an estimate when $a_{kn}$ is supported on $k\leq n$ and given by
$$a_{kn}= \frac{(-1)^k\Gamma(n+2-a)}{k!(n-k)!}\frac{(\lambda+a-\varepsilon_k-k)_n}{(\lambda-\varepsilon_k)_{n+1}}\left(\frac 1{x_1+k}+\frac 1{x_2+k}-\frac 1{\lambda+k}\right),$$
where
$$\varepsilon_k=\frac{(\lambda-x_1)(\lambda-x_2)}{\lambda+k}.$$
Note that $\varepsilon_k\rightarrow 0$ quickly enough so that
\begin{equation}\label{kek}k^{\varepsilon_k}\rightarrow 1,\qquad k\rightarrow\infty.
\end{equation}

It will be useful to write
\begin{align}\nonumber a_{kn}&=\frac{\Gamma(\lambda-\varepsilon_k)}{\Gamma(\lambda+a-\varepsilon_k)\Gamma(1-\lambda-a+\varepsilon_k)}\cdot\frac{\Gamma(n+2-a)\Gamma(\lambda+a-\varepsilon_k+n)}{\Gamma(n+1)\Gamma(\lambda+1-\varepsilon_k+n)} \\
&\nonumber\quad\times\frac{\Gamma(1-\lambda-a+\varepsilon_k+k)}{\Gamma(k+1)}\cdot\frac{(-n)_k}{(1-\lambda-a+\varepsilon_k-n)_k}
 \\
&\label{aknx}\quad\times\left(\frac 1{x_1+k}+\frac 1{x_2+k}-\frac 1{\lambda+k}\right).
 \end{align}
 
 The first factor in \eqref{aknx} has a finite limit as $k\rightarrow\infty$ and is hence bounded by a constant. By \eqref{wsf}, the second factor 
 is bounded by a constant times $n^{1-a}(n-\varepsilon_k)^{a-1}$. However, this is itself a bounded quantity, so the second factor is bounded. By \eqref{gas} and \eqref{kek}, the third factor is 
\begin{equation}\label{tfe}\frac{\Gamma(1-\lambda-a+\varepsilon_k+k)}{\Gamma(k+1)}=\mathcal O(k^{-\lambda-a+\varepsilon_k})=\mathcal O(k^{-\Re(\lambda+a)}). \end{equation}
We write the fourth factor as
\begin{equation}\label{ff}\frac{(-n)_k}{(1-\lambda-a+\varepsilon_k-n)_k}
=\prod_{j=0}^{k-1}\frac{(n-j)}{(\lambda+a-\varepsilon_k-1+n-j)}.\end{equation}
Assuming that $\Re(\lambda+a)>1$, 
we have $\Re(\lambda+a-\varepsilon_k-1)>0$ for large enough 
$k$. Then,
each factor in \eqref{ff} is  bounded in modulus by $1$. 
Hence, \eqref{ff} is bounded by a constant. 
The final factor in  $\eqref{aknx} $ is $\mathcal O(k^{-1})$. This all shows that
$$|a_{kn}|\leq C k^{-1-\Re(\lambda+a)}, $$
with $C$ independent of    $k$ and  $n$. We can hence apply Tannery's theorem and deduce Corollary \ref{cbl} under the assumption $\Re(\lambda+a)>1$.

To weaken the assumption to $\Re(\lambda+a)>0$, we consider the series in \eqref{cdi} as a function of $a$. 
If we can show that it
converges locally uniformly in   $\Re(\lambda+a)>0$,   the general case of Corollary~\ref{cbl} 
follows by analytic continuation. The terms in the relevant series are formed by the first, third and last factor in \eqref{aknx}.
It is clear that the bound on the first factor  can be made locally uniform in $a$. By \eqref{tfe},
 if $\Re(\lambda+a)>\varepsilon>0$, the second factor is $\mathcal O(k^{-\varepsilon})$. Since the final factor is
 $\mathcal O(k^{-1})$, each term can be
 estimated with $C k^{-\varepsilon-1}$, where $C$ is locally bounded as a function of $a$. This completes the proof of Corollary \ref{cbl}.

Let us now turn to Corollary \ref{acc}. Formally, \eqref{clf} is obtained by multiplying \eqref{affp} with $n^{s-2t-u+3}$
and letting $n\rightarrow\infty$. 
In the first sum on the right of \eqref{affp}, the resulting terms can be factored as
\begin{multline}\label{aff}
\frac{\Gamma(\lambda-\varepsilon_k)}{\Gamma(\lambda+u-\varepsilon_k)\Gamma(1-\lambda-u+\varepsilon_k)\Gamma(t-\lambda+\varepsilon_k)}\\
\begin{split}&\times \frac{n^{s-2t-u+3}\Gamma(n+t)\Gamma(\lambda+u-\varepsilon_k+n)\Gamma(t-\lambda+\varepsilon_k+n)}{\Gamma(n+1)\Gamma(\lambda+1-\varepsilon_k+n)\Gamma(s+1-\lambda+\varepsilon_k+n)}\\
&\times\frac{ \Gamma(1-\lambda-u+\varepsilon_k+k)\Gamma(s-\lambda+\varepsilon_k+k)}{\Gamma(k+1)\Gamma(t+k)}\cdot\frac{(-n)_k}{(1-\lambda-u+\varepsilon_k-n)_k}\\
&\times\frac{(t+n)_k}{(s+1-\lambda+\varepsilon_k+n)_k}\cdot\left(\frac1{x_1+k}+\frac1{x_2+k}-\frac1{\lambda+k}\right).
\end{split}\end{multline}
In particular, we have a factor  of the form \eqref{ff} but with $a$ replaced by $u$. Under the  assumption 
$\Re(\lambda+u)>1$, it can be estimated by $1$, for large enough $k$. In order to estimate
$$\frac{(t+n)_k}{(s+1-\lambda+\varepsilon_k+n)_k}
 =\prod_{j=0}^{k-1}\frac{t+n+j}{s+1-\lambda+\varepsilon_k+n+j}, $$
  we will assume that $\Re(\lambda-s+t)>1$.  Note that if $\Re(a)>\Re(b)$, then $|n+a|>|n+b|$ if $n$ is large enough. If we choose $k$ large enough so that $\Re(\lambda-s+t-\varepsilon_k)>1$
and then $n$ large enough so that 
$|t+n+j|>|s+1-\lambda+\varepsilon_k+n+j|$ for all $j$, we have
\begin{multline*}\left|\frac{(t+n)_k}{(s+1-\lambda+\varepsilon_k+n)_k}\right|\leq \left|\frac{(t+n)_n}{(s+1-\lambda+\varepsilon_k+n)_n}\right|\\
=\left|\frac{\Gamma(t+2n)\Gamma(s+1-\lambda+\varepsilon_k+n)}{\Gamma(t+n)\Gamma(s+1-\lambda+\varepsilon_k+2n)}\right|\sim 2^{\Re(\lambda-s+t-\varepsilon_k)-1},\qquad n\rightarrow\infty.
  \end{multline*}
Hence, this factor can be estimated by a constant. The remaining factors in \eqref{aff} can be treated exactly as in \eqref{aknx}.

Turning to the second term on the right-hand side of \eqref{affp}, we note that
$$c_k^\pm=\frac{s+k}2\pm \left(\frac{s+k}2-\lambda+\mathcal O(k^{-1})\right),\qquad k\rightarrow\infty. $$
Hence, 
$$c_k^+=s-\lambda+\delta_k+k,\qquad c_k^-=\lambda-\delta_k, $$
where $k^{\delta_k}\rightarrow 1$ as $k\rightarrow\infty$.
We write the terms as
\begin{multline*}
\frac{\Gamma(\lambda-\delta_k)}{\Gamma(t-\lambda+\delta_k)\Gamma(t-s+\lambda-\delta_k)\Gamma(1+s-t-\lambda+\delta_k)}\\
\begin{split}&\times\frac{\Gamma(s+u+n)\Gamma(t-\lambda+\delta_k+n)\Gamma(t-s+\lambda-\delta_k+n)n^{s-2t-u+3}}{\Gamma(n+1)\Gamma(1+s-\lambda+\delta_k+n)\Gamma(1+\lambda-\delta_k+n)}\\
&\times\frac{\Gamma(s-\lambda+\delta_k+k)\Gamma(1+s-t-\lambda+\delta_k+k)}{\Gamma(k+1)\Gamma(s+u+k)}\\
&\times\frac{(-n)_k}{(1+s-t-\lambda+\delta_k-n)_k}\cdot\frac{(s+u+n)_k}{(1+s-\lambda+\delta_k+n)_k}\cdot \frac1{s+k-x_1-x_2},
\end{split}\end{multline*}
where each factor can be treated as before, under the same assumptions $\Re(\lambda+u)>1$ and
 $\Re(\lambda-s+t)>1$. The analytic continuation to  $\Re(2\lambda-s+t+u)>1$ works as for Corollary \ref{cbl}.

Finally, we turn to Corollary \ref{gcac}. Writing $s=e_1(\mathbf x)$, we have
$$\eta_k^+=s-\lambda+\delta_k+k,\qquad \eta_k^-=\lambda-\delta_k, $$
where $k^{\delta_k}\rightarrow 1$, $k\rightarrow\infty$. We need to multiply \eqref{yl} with $n^{2s+3-3u}$ before letting $n\rightarrow\infty$. The terms on the right can   be factored as
\begin{multline*}
\frac{\Gamma(\lambda-\delta_k)}{\Gamma(u-\lambda+\delta_k)\Gamma(u-s+\lambda-\delta_k)\Gamma(1+s-u-\lambda+\delta_k)}\\
\begin{split}&\times\frac{\Gamma(u+n)\Gamma(u-\lambda+\delta_k+n)\Gamma(u-s+\lambda-\delta_k+n)n^{2s+3-3u}}{\Gamma(n+1)\Gamma(1+s-\lambda+\delta_k+n)\Gamma(1+\lambda-\delta_k+n)}\\
&\times\frac{\Gamma(s-\lambda+\delta_k+k)\Gamma(1+s-u-\lambda+\delta_k+k)}{\Gamma(k+1)\Gamma(u+k)}\\
&\times\frac{(-n)_k}{(1+s-u-\lambda+\delta_k-n)_k}\cdot\frac{(u+n)_k}{(1+s-\lambda+\delta_k+n)_k}\cdot \left(\prod_{j=1}^3\frac1{x_j+k}-\frac 1{s+k}\right).
\end{split}\end{multline*}
We then proceed exactly  as for Corollary \ref{acc}.

 \end{document}